\documentclass[12pt]{iopart}

\usepackage{setstack}
\usepackage{latexsym,amssymb,enumerate,bbm,graphicx,psfrag,color,url}
\usepackage[latin1]{inputenc}

\def\norm#1{\hspace{0.2ex} \|#1\| \hspace{0.2ex}} 

\newcommand{\labeq}[1]{\label{eq:#1}}			
\def\req#1{{\rm(\ref{eq:#1})}}

\newcommand{\R}{\ensuremath{\mathbbm{R}}} 
\newcommand{\Rd}{\ensuremath{\mathbbm{R}_\diamond}} 
\newcommand{\N}{\ensuremath{\mathbbm{N}}}

\newcommand{\dx}[1][x]{\ensuremath{\,{\rm{d}} #1}}

\newcommand{\Ld}{L_\diamond}
\newcommand{\Hd}{H_\diamond}
\newcommand{\LL}{\mathcal L}

\newcommand{\FF}{\mathcal F}
\newcommand{\KK}{\mathcal K}
\newcommand{\Langle}{\left\langle}
\newcommand{\Rangle}{\right\rangle}

\newcommand{\out}{\mathrm{out}_{\Sigma}\,}
\newcommand{\inn}{\mathrm{int}\,}
\renewcommand{\span}{\mathrm{span}\,}

\newcommand{\kommentar}[1]{}

\newtheorem{theorem}{Theorem}[section]
\newtheorem{lemma}[theorem]{Lemma}

\newtheorem{definition}[theorem]{Definition}
\newtheorem{remark}[theorem]{Remark}

\begin{document}
\title[Uniqueness and Lipschitz stability in EIT with finitely many electrodes]{Uniqueness and Lipschitz stability in Electrical Impedance Tomography with finitely many electrodes}

\renewcommand{\footnoterule}{%
  \kern -3pt
  \hrule width \textwidth height 1pt
  \kern 2pt
}
\footnotetext{%
\scriptsize
This is the version of the article before peer review or editing, as submitted by an
author to \emph{Inverse Problems}. IOP Publishing Ltd is not responsible for any
errors or omissions in this version of the manuscript or any version derived from
it. The Version of Record is available online at \url{https://doi.org/10.1088/1361-6420/aaf6fc}.
}     

\author{Bastian Harrach$^1$}
\address{$^1$ Institute for Mathematics, Goethe-University Frankfurt, Frankfurt am Main, 
Germany}
\ead{harrach@math.uni-frankfurt.de}

\begin{abstract}
For the linearized reconstruction problem in Electrical Impedance Tomography (EIT) with the Complete Electrode Model (CEM), Lechleiter and Rieder (2008 Inverse Problems 24 065009) have shown that a piecewise polynomial conductivity on a fixed partition is uniquely determined if enough electrodes are being used. We extend their result to the full non-linear case and show that measurements on a sufficiently high number of 
electrodes uniquely determine a conductivity in any finite-dimensional subset of piecewise-analytic functions. We also prove Lipschitz stability, and derive analogue results for the continuum model, where finitely many measurements determine a finite-dimensional Galerkin projection
of the Neumann-to-Dirichlet operator on a boundary part. 
\end{abstract}




\section{Introduction}

We consider the inverse conductivity problem of determining the coefficient function $\sigma$ in the elliptic partial differential equation
\begin{equation}\labeq{intro:EIT}
\nabla \cdot (\sigma \nabla u) =0 \quad \mbox{ in } \Omega
\end{equation}
from knowledge of boundary measurements of $u$. The problem arises in Electrical Impedance Tomography (EIT), or electrical resistivity tomography,  
which is a novel technique to image the conductivity distribution $\sigma$ inside a subject $\Omega$ from electric voltage and current measurements on the subject's boundary $\partial \Omega$, cf.\ \cite{henderson1978impedance,barber1984applied,wexler1985impedance,newell1988electric,metherall1996three,cheney1999electrical,borcea2002electrical,borcea2003addendum,lionheart2003eit,holder2004electrical,bayford2006bioimpedance,uhlmann2009electrical,adler2015electrical,martinsen2011bioimpedance,seo2013electrical}, and the references therein for a broad overview on the developments in EIT.

To model the boundary measurements we consider the continuum model, where we measure the local Neumann-to-Dirichlet operator (on a boundary part $\Sigma\subseteq \partial \Omega$)
\[
\Lambda(\sigma):\ g\mapsto u|_{\Sigma},\quad \mbox{ where $u$ solves \req{intro:EIT} with }
\sigma\partial_\nu u|_{\partial \Omega}=\left\{ \begin{array}{l l} g & \mbox{ on $\Sigma$,}\\ 0 & \mbox{ else,}\end{array}\right.
\] and the more realistic Complete Electrode Model (CEM) with electrodes $E_1,\ldots,E_M\subseteq \partial \Omega$ all having the same contact impedance $z>0$.
In the CEM, we measure
\begin{eqnarray*}
R_M(\sigma):\ (J_1,\ldots,J_M)\mapsto (U_1,\ldots,U_M),
\end{eqnarray*}
where $u$ solves \req{intro:EIT} with 
\begin{eqnarray*}
\sigma\partial_\nu u=0 \quad & \mbox{on $\partial \Omega\setminus \bigcup_{m=1}^M E_m$,}\\
u+z \sigma \partial_\nu u=\mathrm{const.}=:U_m \quad & \mbox{on $E_m$, \quad $m=1,\ldots,M$,}\\
\int_{E_m} \sigma \partial_\nu u|_{E_m}\dx[s] =J_m \quad & \mbox{on $E_m$, \quad $m=1,\ldots,M$.}
\end{eqnarray*}

The question whether full or local Neumann-Dirichlet-measurements uniquely determine the coefficient function $\sigma$ has become famous under the name Calder\'on problem \cite{calderon1980inverse,calderon2006inverse}, and has been intensively studied in the mathematical literature due to its practical relevance for EIT and many other related inverse coefficient problems, cf.\
\cite{kohn1984determining,kohn1985determining,druskin1998uniqueness,sylvester1987global,nachman1996global,ammari2004reconstruction,astala2006calderon,isakov2007uniqueness,kenig2007calderon,imanuvilov2010calderon,haberman2013uniqueness,guillarmou2013calderon,kenig2014calderon,kenig2014recent,imanuvilov2015neumann,caro2016global,krupchyk2016calderon}.

In this work we will study the question whether $\sigma$ can be uniquely and stably reconstructed from a finite number of electrode measurements.
A natural discretization is to assume that $\sigma$ is piecewise constant (or piecewise polynomial) on a given resolution or partition of $\Omega$,
so that $\sigma$ will lie in an a-priori known finite-dimensional subset $\FF$ of piecewise-analytic functions. Moreover, it seems natural to
assume that upper and lower bounds on the conductivity are a-priori known, i.e.,
\[
\sigma\in \FF_{[a,b]}:=\{ \sigma\in \FF:\ a\leq  \sigma(x)\leq  b\mbox{ for all } x\in \Omega\}.
\]
Our main result for the continuum model is that a (sufficiently high dimensional) finite-dimensional Galerkin projection $G_N \Lambda(\sigma) G_N^*$ already uniquely determines $\sigma$ and that Lipschitz stability holds
\[
\exists c>0:\ c \norm{\sigma_1-\sigma_2}\leq \norm{ G_N \left(\Lambda(\sigma_1)-\Lambda(\sigma_2)\right) G_N^*},
\]
cf.\ Theorem~\ref{thm:stability_Galerkin}.

Under the additional assumption that $\sigma$ is an a-priori known smooth function close to the boundary, we then 
turn to the Complete Electrode Model. We show that a (sufficiently large) finite number of electrodes suffices to uniquely 
determine $\sigma$ with Lipschitz stability
\[
\exists c>0:\ c \norm{\sigma_1-\sigma_2}\leq \norm{ R_M(\sigma_1)-R_M(\sigma_2)},
\]
cf.\ Theorem~\ref{thm:stability_CEM}. This shows that the discretized EIT problem is uniquely and stably solvable if enough electrodes are being used,
which may be relevant for practical implementations of EIT reconstruction algorithms.

Note that our results are non-constructive, we do not have a practically useful estimate of
the Lipschitz constant or the required number of electrodes yet. Also note, that the necessary number of electrodes and the stability constant $c>0$ depend on the ansatz set $\FF_{[a,b]}$. Due to the intrinsic ill-posedness of the non-discretized EIT problem, we can naturally expect that a larger set $\FF_{[a,b]}$ will lead to worse stability constants
and a higher required number of electrodes, with $c\to 0$ and $M\to \infty$ when $\mathrm{dim}(\span \FF)\to \infty$. 

Let us give some more references on related results and the origins of our approach. A recent preprint of Alberti and Santacesaria \cite{alberti2018calder}
uses complex geometrical optics solutions to show that (in the continuuum model) there exists a finite number of boundary voltages, so that the knowledge of the corresponding boundary currents uniquely determines the conductivity $\sigma$ and that Lipschitz stability holds. Their result holds in dimension $d\geq 3$ with measurements on the full boundary $\partial \Omega$, $\sigma$ is assumed to be identically one close to $\partial \Omega$, bounded by a-priori known constants, and $\frac{\Delta \sigma}{\sigma}$ has to belong to an a-priori known finite-dimensional subspace of $L^\infty$. Our result in this work works with less restrictive assumptions
as we can treat any dimension $d\geq 2$, partial boundary data, and the complete electrode model. But, on the other hand, we require the 
assumption of piecewise-analyticity which is more restrictive than the assumptions in \cite{alberti2018calder}. 

For the linearized EIT problem (both, in the continuum model, and with the CEM), Lechleiter and Rieder \cite{lechleiter2008newton} have shown that a piecewise polynomial conductivity on a fixed partition is uniquely determined if enough electrodes are being used. The main tool in \cite{lechleiter2008newton} is the theory of localized potentials devoloped by the author \cite{gebauer2008localized} and the 
convergence of CEM-solutions to solutions of the continuum model shown by Hyv\"onen, Lechleiter and Hakula \cite{hyvonen2004complete,lechleiter2008factorization}. Our result uses similar tools and first treats the non-linear EIT problem with the continuum model using localized potentials \cite{gebauer2008localized,harrach2013monotonicity} and monotonicity estimates between the non-linearized and the linearized problem from Ikehata, Kang, Seo and Sheen \cite{kang1997inverse,ikehata1998size}. Then we extend the results to the CEM using recent results on the approximation of the continuum model by the CEM from Hyvönen, Garde and Staboulis \cite{hyvonen2009approximating,garde2017convergence}. 

The idea of using monotonicity estimates and localized potentials techniques has lead to a number of results for inverse coefficient problems \cite{harrach2009uniqueness,harrach2010exact,harrach2012simultaneous,arnold2013unique,harrach2013monotonicity,barth2017detecting,harrach2017local,brander2018monotonicity,griesmaier2018monotonicity,harrach2018helmholtz,harrach2018fractional,harrach2018localizing},
and several recent works build practical reconstruction methods on monotonicity properties
\cite{tamburrino2002new,harrach2015combining,harrach2015resolution,harrach2016enhancing,maffucci2016novel,tamburrino2016monotonicity,garde2017comparison,garde2017convergence,garde2017regularized,su2017monotonicity,ventre2017design,harrach2018monotonicity,zhou2018monotonicity}.
Together with the recent preprint \cite{harrach2018global}, the present work shows that this idea 
can also be used to obtain Lipschitz stability estimates, which are usually derived from technically more challenging approaches involving Carleman estimates or quantitative unique continuation, cf.\ \cite{kazemi1993stability,alessandrini1996determining,imanuvilov1998lipschitz,imanuvilov2001global,cheng2003lipschitz,alessandrini2005lipschitz,bacchelli2006lipschitz,bellassoued2006lipschitz,klibanov2006lipschitz,bellassoued2007lipschitz,klibanov2006lipschitz,klibanov2006lipschitz_nonstandard,sincich2007lipschitz,yuan2007lipschitz,yuan2009lipschitz,beretta2011lipschitz,beretta2013lipschitz,melendez2013lipschitz,alessandrini2017lipschitz,beretta2017uniqueness,beilina2017lipschitz,alessandrini2018lipschitz,ruland2018lipschitz}.

The work is organized as follows. In section \ref{section:continuum} we treat the continuum model,
and show that the Neumann-to-Dirichlet operator or a (sufficiently high dimensional) finite-dimensional Galerkin projection 
uniquely determine the conductivity with Lipschitz stability. We formulate our main results for the continuum model in Theorem~\ref{thm:stability_full} and Theorem \ref{thm:stability_Galerkin} in subsection \ref{subsect:continuum_setting}, 
summarize some known results from the literature in subsection \ref{subsect:continuum_tools}, and the prove the theorems in subsection~\ref{subsect:proofs}. In section \ref{section:CEM} we then treat the Complete Electrode Model. Again we first formulate a uniqueness and Lipschitz stability result in Theorem~\ref{thm:stability_CEM} in subsection~\ref{subsect:setting}, then summarize known results from the literature in subsection \ref{subsect:CEM_tools}, and finally prove the theorem in subsection~\ref{subsect:proofs_CEM}

\section{Uniqueness and Lipschitz stability from continuous data}\label{section:continuum}

\subsection{Setting and main results}\label{subsect:continuum_setting}

Let $\Omega\subset\R^d$, $d\geq 2$ be a bounded domain with smooth boundary $\partial\Omega$ and outer normal vector $\nu$. $L_+^\infty(\Omega)$ denotes the subspace of $L^\infty(\Omega)$-functions with positive essential infima. $\Hd^1(\Omega)$ and $\Ld^2(\partial \Omega)$ denote the spaces of $H^1$- and $L^2$-functions with vanishing integral mean on $\partial \Omega$.

For $\sigma\in L_+^\infty(\Omega)$, and a relatively open boundary part $\Sigma\subseteq \partial \Omega$, 
the local Neumann-to-Dirichlet (NtD) operator $\Lambda(\sigma)$ is defined by 
\[
\Lambda(\sigma):\ L^2_\diamond(\Sigma)\to L^2_\diamond(\Sigma), \quad g\mapsto u^g_\sigma|_{\Sigma},
\]
where $u^g_\sigma\in H^1_\diamond(\Omega)$ is the unique solution of
\begin{equation}\labeq{math_model}
\nabla\cdot\left(\sigma\nabla u^g_\sigma\right) = 0 \mbox{ in }\Omega,\quad\sigma\partial_\nu u_\sigma^g\vert_{\partial\Omega} = \left\{ \begin{array}{l l} g & \mbox{ on }\Sigma,\\
0 & \mbox{ else.}\end{array}\right.
\end{equation}
This is equivalent to the variational formulation that $u^g_\sigma\in H^1_\diamond(\Omega)$ solves
\begin{equation}\labeq{EIT_CM_variational}
\int_\Omega \sigma \nabla u^g_\sigma \cdot \nabla w \dx = \int_{\Sigma} g w|_{\Sigma}\dx[s] \quad \mbox{ for all } w\in \Hd^1(\Omega).
\end{equation}
It is well known and easily shown that $\Lambda(\sigma)$ is compact and self-adjoint. 

We will consider conductivities that are a-priori known to belong to a finite dimensional set of piecewise-analytic 
functions and that are bounded from above and below by a-priori known constants. To that end, we first define piecewise-analyticity as in 
\cite[Def.~2.1]{harrach2013monotonicity}:
\begin{definition}
\label{def:piecewise}
\begin{enumerate}[(a)]
\item A subset $\Gamma\subseteq \partial O$ of the boundary of an open set $O\subseteq \R^n$ is called a \emph{smooth boundary piece} if
it is a $C^\infty$-surface and $O$ lies on one side of it, i.e., if for each $z\in \Gamma$
there exists a ball $B_\epsilon(z)$ and a function $\gamma\in C^\infty(\R^{n-1},\R)$ such that upon relabeling and reorienting
\begin{eqnarray*}
\Gamma=\partial O\cap B_\epsilon(z)=\{ x\in B_\epsilon(z) \; | \; x_n=\gamma(x_1,\ldots,x_{n-1}) \},\\
O\cap B_\epsilon(z)=\{ x\in B_\epsilon(z) \; | \; x_n>\gamma(x_1,\ldots,x_{n-1}) \}.
\end{eqnarray*}
\item $O$ is said to have \emph{smooth boundary} if $\partial O$ is a union of smooth boundary pieces.
$O$ is said to have \emph{piecewise smooth boundary} if $\partial O$ is a countable union of the closures of smooth boundary pieces.
\item A function $\kappa \in L^\infty(\Omega)$ is called \emph{piecewise analytic} if there exist
finitely many pairwise disjoint subdomains $O_1,\ldots,O_M\subset \Omega$ with piecewise
smooth boundaries, such that $\overline{\Omega}= \overline{O_1\cup \ldots \cup O_M}$,
and $\kappa|_{O_m}$ has an extension which is (real-)analytic in a neighborhood of $\overline{O_m}$,
$m=1,\ldots,M$.
\end{enumerate}
\end{definition}

Note that (to the knowledge of the author), it is not clear whether the sum of two piecewise-analytic functions is always piece\-wise-analytic, 
i.e., whether the set of piecewise-analytic functions is a vector space. But finite-dimensional vector spaces of piecewise-analytic functions (or subsets thereof) naturally arise as parameter spaces for the inverse conductivity problem, e.g., when we fix a partition of the imaging domain $\Omega$ into a finite number of subdomains (e.g., triangles, pixels, or voxels) and the conductivity is assumed to be a polynomial of fixed maximal order on each of these subdomains. Therefore, we make the following definition:
\begin{definition}
\label{def:findim_subsets}
A set $\FF\subseteq L^\infty(\Omega)$ is called a \emph{finite-di\-men\-sional subset of piecewise-analytic functions} if its linear span
\[
\span \FF=\left\{ \sum_{j=1}^k \lambda_j f_j:\ k\in \N,\ \lambda_j\in \R,\ f_j\in \FF\right\}\subseteq L^\infty(\Omega)
\] 
contains only piecewise-analytic functions and $\mathrm{dim}(\span\FF)<\infty$.
\end{definition}

Given a finite-dimensional subset $\FF$ of piecewise analytic functions and
two numbers $b>a>0$, we denote the set 
\[
\FF_{[a,b]}:=\{ \sigma\in \FF:\ a\leq  \sigma(x)\leq  b\mbox{ for all } x\in \Omega\}.
\]

Throughout this paper, the domain $\Omega$, the finite-dimensional subset $\FF$ and the bounds $b>a>0$ are fixed, and the constants in the Lipschitz stability results will depend on them.


Our first result shows Lipschitz stability for the inverse conductivity problem in $\FF_{[a,b]}$ 
when the complete infinite-dimensional Neumann-to-Dirichlet-operator is measured.

\begin{theorem}\label{thm:stability_full}
There exists $c>0$ such that
\[
\norm{\Lambda(\sigma_1)-\Lambda(\sigma_2)}_{\LL(\Ld^2(\Sigma))}\geq c \norm{\sigma_1-\sigma_2}_{L^\infty(\Omega)} \quad \mbox{ for all } \sigma_1,\sigma_2\in \FF_{[a,b]}.
\]
\end{theorem}
\paragraph{Proof.}
Theorem \ref{thm:stability_full} will be proven in subsection \ref{subsect:proofs}.
\hfill $\Box$

We then turn to the question whether $\sigma\in \FF_{[a,b]}$ is already uniquely determined by finitely many boundary measurements in the continuum model.
For a (finite- or infinite-dimensional) subspace $G\subseteq \Ld^2(\Sigma)$ we denote by
\[
P_G:\ \Ld^2(\Sigma) \to G, \quad P_G g=\left\{ \begin{array}{l l} g & \mbox{ if } g\in G,\\ 0 & \mbox{ if } g\in G^\perp.\end{array} \right.
\]
the orthogonal projection operator on $G$ with respect to the $L^2$-scalar product
\begin{equation}\labeq{L2_scalar_product}
\langle g,h\rangle:=\int_{\Sigma} gh \dx[s] \quad \mbox{ for all } g,h\in L^2(\Sigma).
\end{equation}
If $G$ is finite dimensional with a basis $G=\span(g_1,\ldots,g_n)$ then measurements of
\[
\langle g_j,\Lambda(\sigma)g_k\rangle \quad j,k=1,\ldots,n
\]
 determine the Galerkin projection of the Neumann-to-Dirichlet operator $P_{G} \Lambda(\sigma) P_G^*$,
so that this can be regarded as a model for finitely many voltage/current measurements in the continuum model.

Our next result shows that this uniquely determines $\sigma\in \FF_{[a,b]}$ (with Lipschitz stability) if the space $G$ is large enough.

\begin{theorem}\label{thm:stability_Galerkin}
For each sequence of subspaces
\[
G_1\subseteq G_2\subseteq\ldots \Ld^2(\Sigma)\quad \mbox{ with } \quad 
\overline{\bigcup_{n\in \N} G_n}=\Ld^2(\Sigma)
\]
there exists $N\in \N$, and $c>0$ such that
\[
\norm{P_{G_n} \left(\Lambda(\sigma_1)-\Lambda(\sigma_2)\right) P_{G_n}^*}_{\LL(\Ld^2(\Sigma))}\geq c \norm{\sigma_1-\sigma_2}_{L^\infty(\Omega)}  
\]
for all $\sigma_1,\sigma_2\in \FF_{[a,b]}$, and all $n\geq N$. 

In particular, this implies that for all $\sigma_1,\sigma_2\in \FF_{[a,b]}$ and all $n\geq N$
\[
P_{G_n} \Lambda(\sigma_1) P_{G_n}^*=P_{G_n} \Lambda(\sigma_2)P_{G_n}^* \quad \mbox{ if and only if } \quad \sigma_1=\sigma_2.
\]
\end{theorem}
\paragraph{Proof.} 
Theorem \ref{thm:stability_Galerkin} will be proven in subsection \ref{subsect:proofs}.
\hfill $\Box$

\subsection{Differentiability, monotonicity and localized potentials}\label{subsect:continuum_tools}

In this subsection, we summarize some known results from the literature, that we will use to prove Theorem \ref{thm:stability_full} and \ref{thm:stability_Galerkin}. As defined in \req{L2_scalar_product},  $\langle \cdot,\cdot \rangle$ always denotes the $L^2(\Sigma)$-scalar product, and
$u^g_\sigma\in H^1_\diamond(\Omega)$ denotes the solution of \req{math_model} with conductivity $\sigma\in L^\infty_+(\Omega)$ and Neumann data $g\in \Ld^2(\Sigma)$ in the following. 

Our first tool is that the Neumann-to-Dirichlet (NtD) operator is continuously Fr\'echet differentiable with respect to the conductivity.
\begin{lemma}\label{lemma:Frechet_differentiable}
\begin{enumerate}[(a)]
\item The mapping 
\[
\Lambda:\ L^\infty_+(\Omega)\to \LL(\Ld^2(\Sigma)),\quad \sigma\mapsto \Lambda(\sigma)
\]
is Fr\'echet differentiable. Its derivative is given by 
\begin{equation}\labeq{Diff_CM}
\Lambda'(\sigma)\in \LL(L^\infty(\Omega), \LL(\Ld^2(\Sigma))), \quad \left(\Lambda'(\sigma)\kappa\right)g=v|_{\Sigma},
\end{equation}
where $v\in \Hd^1(\Omega)$ solves
\[
\int_\Omega \sigma \nabla v\cdot \nabla w \dx = - \int_\Omega \kappa \nabla u_\sigma^{g} \cdot \nabla w \dx \quad \mbox{ for all } w\in \Hd^1(\Omega).
\]
\item For all $\sigma\in L^\infty_+(\Omega)$ and $\kappa\in L^\infty(\Omega)$ the operator $\Lambda'(\sigma)\kappa\in  \LL(\Ld^2(\Sigma))$ is self-adjoint and compact,
and it fulfills
\[
\langle \left(\Lambda'(\sigma)\kappa\right) g,h\rangle=\int_\Omega \sigma \nabla u_\sigma^h \cdot \nabla u_\sigma^g \dx = -\int_\Omega \kappa \nabla u_\sigma^g \cdot \nabla u_\sigma^h \dx.
\]
for all $g,h\in \Ld^2(\Sigma)$.
\item The mapping
\[
\Lambda':\ L^\infty_+(\Omega)\to \LL(L^\infty(\Omega), \LL(\Ld^2(\Sigma))),\quad \sigma\mapsto \Lambda'(\sigma)
\]
is continuous.
\end{enumerate}
\end{lemma}
\paragraph{Proof.}
This follows from the variational formulation of the conductivity equation \req{EIT_CM_variational}, cf., e.g., \cite[Section~2]{lechleiter2008newton} or \cite[Appendix~B]{garde2017convergence}.
\hfill $\Box$

Our next tool is a monotonicity relation between the NtD-operator and its derivative that goes back to Ikehata, Kang, Seo, and Sheen \cite{kang1997inverse,ikehata1998size}, and has been used in several other works, cf.\ the list of works on monotonicity-based methods
cited in the introduction.
\begin{lemma}\label{lemma:monotony}
For all $\sigma_1,\sigma_2\in L^\infty_+(\Omega)$ and $g\in L^2_\diamond(\Sigma)$, it holds that
\begin{eqnarray}
\nonumber \langle \left(\Lambda'(\sigma_2)(\sigma_1-\sigma_2)\right) g,g\rangle
&=\int_\Omega(\sigma_2-\sigma_1) |\nabla u_{\sigma_2}^g|^2 \dx\\
\label{eq:monotony_inequality}
&\leq \Langle g, \left(\Lambda(\sigma_1)-\Lambda(\sigma_2)\right)g \Rangle.
\end{eqnarray}
\end{lemma}
\paragraph{Proof.}
See, e.g., \cite[lemma~2.1]{harrach2010exact}.
\hfill $\Box$

The energy terms $|\nabla u_\sigma^g|^2$ in the monotonicity estimate can be controlled using the technique of localized potentials \cite{gebauer2008localized}. Roughly speaking, the energy $|\nabla u_\sigma^g|^2$ can be made arbitrarily 
large in a subset $D_1\subseteq \Omega$ without making it large in another subset $D_2\subseteq \Omega$ whenever $D_1$ can be reached from
the boundary without passing $D_2$. 

To formulate this rigorously, we adopt the notation from \cite[Def.~2.2, 2.3]{harrach2013monotonicity} and denote by 
$\inn D$ the topological interior of a subset $D\subseteq \overline\Omega$, and by $\out D$ its outer hull, i.e. 
\begin{eqnarray*}
\out D:=\overline \Omega \setminus \bigcup \left\{ 
U\subseteq \overline \Omega:\  \mbox{$U$ rel.\ open,\ $U\cap \Omega$ connected,
\ $U\cap \Sigma\neq \emptyset$}
\right\}.
\end{eqnarray*}
With this notation, we have the following localized potentials result:

\begin{lemma}\label{lemma:loc_pot}
Let $\sigma\in L^\infty_+(\Omega)$ be piecewise analytic and let $D_1,D_2\subseteq\overline\Omega$ be two measurable sets with 
\[
\inn D_1 \nsubseteq \out D_2.
\]

Then there exists a sequence of currents $(g_m)_{m\in\mathbb{N}}\subset L^2_\diamond(\Sigma)$ such that the corresponding
solutions $(u_\sigma^{g_m})_{m\in \mathbb{N}}\subset \Hd^1(\Omega)$ 
fulfill
\[
 \lim_{m\rightarrow\infty}\int_{D_1}\left|\nabla u_\sigma^{g_m}\right|^2\,\mathrm{d}x=\infty\quad\mbox{and}\quad\lim_{m\rightarrow\infty}\int_{D_2}\left|\nabla u_\sigma^{g_m}\right|^2\,\mathrm{d}x=0.
\]
\end{lemma}
\paragraph{Proof.}
\cite[Thm.~3.6 and Sect.~4.3]{harrach2013monotonicity}
\hfill $\Box$

We will also need the following definiteness property of piecewise-analytic functions from \cite{harrach2013monotonicity}:
\begin{lemma}\label{lemma:pcw_anal_definiteness}
Let $0\not\equiv\kappa\in L^\infty(\Omega)$ be piecewise-analytic. Then there exist two sets 
\[
D_1=\inn D_1\quad \mbox{ and } \quad D_2=\out D_2
\]
(i.e., $D_1$ is open, $D_2$ is closed, $\Omega\setminus D_2$ is connected, and $\Sigma \cap \overline \Omega\setminus D_2\neq \emptyset$)
with
\[
D_1=\inn D_1 \nsubseteq \out D_2=D_2
\]
and either 
\begin{enumerate}[(i)]
\item $\kappa|_{\Omega\setminus D_2}\geq 0$ and $\kappa|_{D_1}\in L_+^\infty(D_1)$, or
\item $\kappa|_{\Omega\setminus D_2}\leq 0$ and $-\kappa|_{D_1}\in L_+^\infty(D_1)$.
\end{enumerate}
\end{lemma}
\paragraph{Proof.}
\cite[Thm.~A.1, Cor.~A.2, and Sect.~4.3]{harrach2013monotonicity}
\hfill $\Box$

\subsection{Proof of Theorem \ref{thm:stability_full} and Theorem \ref{thm:stability_Galerkin}}\label{subsect:proofs}

We can now prove Theorem \ref{thm:stability_full} and Theorem \ref{thm:stability_Galerkin}. For the sake of brevity, we write
$\norm{\cdot}$ for $\norm{\cdot}_{\LL(\Ld^2(\Sigma))}$, $\norm{\cdot}_{L^\infty(\Omega)}$ and
$\norm{\cdot}_{\Ld^2(\Sigma)}$ throughout this subsection.

We follow the approach in \cite{harrach2018global} and first use the monotonicity relation in lemma~\ref{lemma:monotony} to bound the difference of the non-linear Neumann-to-Dirichlet operators by an expression containing their linearized counterparts.

\begin{lemma}\label{lemma:estimate_with_f}
For all $\sigma_1,\sigma_2\in \FF_{[a,b]}$ with $\sigma_1\not\equiv\sigma_2$, 
\[
\frac{\norm{\Lambda(\sigma_1)-\Lambda(\sigma_2)}}
{\norm{\sigma_1-\sigma_2}}
\geq   \inf_{(\tau_1,\tau_2,\kappa)\atop \in \FF_{[a,b]}\times \FF_{[a,b]}\times \KK} \ \sup_{g\in \Ld^2(\Sigma), \norm{g}=1} f(\tau_1,\tau_2,\kappa,g),
\]
where $f:\ L^\infty_+(\Omega)\times L^\infty_+(\Omega)\times L^\infty(\Omega)\times \Ld^2(\Sigma)\to \R$ is defined by
\[
f(\tau_1,\tau_2,\kappa,g):=\max \left\{ \Langle \left(\Lambda'(\tau_1) \kappa\right) g,g \Rangle,
-\Langle \left(\Lambda'(\tau_2) \kappa\right) g,g \Rangle\right\},
\]
and $\KK:=\{ \kappa\in \span\FF:\ \norm{\kappa}=1\}$.
\end{lemma}
\paragraph{Proof.}
The Neumann-to-Dirichlet-operators are self-adjoint so that for all $\sigma_1,\sigma_2\in L^\infty(\Omega)$
\[
\norm{\Lambda(\sigma_1)-\Lambda(\sigma_2)}=\sup_{g\in \Ld^2(\Sigma), \norm{g}=1} \left| \Langle g, \left(\Lambda(\sigma_1)-\Lambda(\sigma_2)\right)g \Rangle \right|.
\]
Using the monotonicity inequality in Lemma~\ref{lemma:monotony} also with interchanged roles of $\sigma_1$ and $\sigma_2$, we obtain that
for all $\sigma_1,\sigma_2\in L^\infty_+(\Omega)$, $\sigma_1\not\equiv \sigma_2$, and all $g\in \Ld^2(\partial \Omega)$
\begin{eqnarray*}
\lefteqn{\left| \Langle g, \left(\Lambda(\sigma_1)-\Lambda(\sigma_2)\right)g \Rangle \right|}\\
&= \max \{ \Langle g, \left(\Lambda(\sigma_2)-\Lambda(\sigma_1)\right)g \Rangle, \Langle g, \left(\Lambda(\sigma_1)-\Lambda(\sigma_2)\right)g \Rangle\}\\
&\geq \max \{ \Langle \Lambda'(\sigma_1) (\sigma_2-\sigma_1) g,g \Rangle,
\Langle \Lambda'(\sigma_2) (\sigma_1-\sigma_2) g,g \Rangle\}\\
&= \norm{\sigma_1-\sigma_2}  \max \left\{ \Langle \Lambda'(\sigma_1) \frac{\sigma_2-\sigma_1}{\norm{\sigma_1-\sigma_2}} g,g \Rangle,
\Langle \Lambda'(\sigma_2) \frac{\sigma_1-\sigma_2}{\norm{\sigma_1-\sigma_2}} g,g \Rangle\right\}\\
&=  \norm{\sigma_1-\sigma_2} f(\sigma_1,\sigma_2,\frac{\sigma_2-\sigma_1}{\norm{\sigma_1-\sigma_2}},g).
\end{eqnarray*}
Hence,
\begin{eqnarray*}
\frac{\norm{\Lambda(\sigma_1)-\Lambda(\sigma_2)}}
{\norm{\sigma_1-\sigma_2}}
&= \sup_{g\in \Ld^2(\Sigma), \norm{g}=1} \frac{\left| \Langle g, \left(\Lambda(\sigma_1)-\Lambda(\sigma_2)\right)g \Rangle \right|}
{\norm{\sigma_1-\sigma_2}}\\
&\geq \sup_{g\in \Ld^2(\Sigma), \norm{g}=1} f(\sigma_1,\sigma_2,\frac{\sigma_2-\sigma_1}{\norm{\sigma_1-\sigma_2}},g)\\
&\geq \inf_{(\tau_1,\tau_2,\kappa)\atop \in \FF_{[a,b]}\times \FF_{[a,b]}\times \KK}\
\sup_{g\in \Ld^2(\Sigma), \norm{g}=1} f(\tau_1,\tau_2,\kappa,g).
\end{eqnarray*}
\hfill $\Box$

Now we use a compactness argument to show that the expression in the lower bound in lemma~\ref{lemma:estimate_with_f} attains its minimum.

\begin{lemma}\label{CM_inf_is_min}
There exists $(\hat \tau_1,\hat \tau_2,\hat \kappa) \in \FF_{[a,b]}\times \FF_{[a,b]}\times \KK$ so that
\[
\inf_{(\tau_1,\tau_2,\kappa)\atop \in \FF_{[a,b]}\times \FF_{[a,b]}\times \KK} \ \sup_{g\in \Ld^2(\Sigma), \norm{g}=1} f(\tau_1,\tau_2,\kappa,g)
= \sup_{g\in \Ld^2(\Sigma), \norm{g}=1} f(\hat \tau_1,\hat \tau_2,\hat \kappa,g).
\]
\end{lemma}
\paragraph{Proof.}
Since $f$ is continuous by lemma \ref{lemma:Frechet_differentiable}, 
the function
\[
(\tau_1,\tau_2,\kappa)\mapsto \sup_{g\in \Ld^2(\Sigma), \norm{g}=1} f(\tau_1,\tau_2,\kappa,g)
\]
is lower semicontinuous and thus attains its minimum over the compact set $ \FF_{[a,b]}\times \FF_{[a,b]}\times \KK$.
\hfill $\Box$

It remains to show that the minimum attained in lemma~\ref{CM_inf_is_min} must be positive. To show that we use the localized potentials from lemma \ref{lemma:loc_pot}.

\begin{lemma}\label{lemma:linearized_uniform_sign}
Let $0\not\equiv\kappa\in L^\infty(\Omega)$ be piecewise-analytic. Then at least one of the following two properties holds true:
\begin{enumerate}[(i)]
\item For all piecewise analytic $\sigma\in L^\infty_+(\Omega)$ there exists $g\in \Ld^2(\partial \Omega)$ with
\[
 -\langle \left(\Lambda'(\sigma)\kappa \right) g,g\rangle
=\int_\Omega \kappa |\nabla u_{\sigma}^g|^2 \dx >0.
\]
\item For all piecewise analytic $\sigma\in L^\infty_+(\Omega)$ there exists $g\in \Ld^2(\partial \Omega)$ with
\[
 -\langle \left(\Lambda'(\sigma)\kappa \right) g,g\rangle
=\int_\Omega \kappa |\nabla u_{\sigma}^g|^2 \dx <0.
\]
\end{enumerate}
Hence, a fortiori, 
\[
\sup_{g\in \Ld^2(\Sigma), \norm{g}=1} f(\tau_1,\tau_2,\kappa,g)>0 \quad \mbox{ for all } (\tau_1,\tau_2,\kappa)\in \FF_{[a,b]}\times \FF_{[a,b]}\times \KK.
\]
\end{lemma}
\paragraph{Proof.}
Using the definiteness property of piecewise analytic functions from lemma~\ref{lemma:pcw_anal_definiteness} we obtain two sets $D_1=\inn D_1$ and $D_2=\out D_2$ with
\[
D_1=\inn D_1 \nsubseteq \out D_2=D_2
\]
and either 
\begin{enumerate}[(i)]
\item $\kappa|_{\Omega\setminus D_2}\geq 0$ and $\kappa|_{D_1}\in L_+^\infty(D_1)$, or
\item $\kappa|_{\Omega\setminus D_2}\leq 0$ and $-\kappa|_{D_1}\in L_+^\infty(D_1)$.
\end{enumerate}
Let $(g_m)_{m\in\mathbb{N}}\subset L^2_\diamond(\Sigma)$ be the localized potentials sequence from 
lemma~\ref{lemma:loc_pot}. Then, in case (a), we obtain 
\begin{eqnarray*}
\lefteqn{-\langle \left(\Lambda'(\sigma)\kappa \right) g_m,g_m\rangle=\int_\Omega \kappa |\nabla u_\sigma^{g_m}|^2 \dx}\\
&= \int_{D_1} \kappa |\nabla u_{\sigma}^{g_m}|^2 \dx  
+ \int_{D_2} \kappa |\nabla u_{\sigma}^{g_m}|^2 \dx 
\int_{\Omega\setminus (D_1\cup D_2)} \kappa |\nabla u_{\sigma}^{g_m}|^2 \dx\\
&\geq \mathrm{ess\,inf}\, \kappa|_{D_1} \int_{D_1} |\nabla u_{\sigma}^{g_m}|^2 \dx
- \norm{\kappa}_{L^\infty(D_2)} \int_{D_2} |\nabla u_{\sigma}^{g_m}|^2 \dx\to \infty,
\end{eqnarray*}
so that $\int_\Omega \kappa |\nabla u_\sigma^{g_m}|^2 \dx>0$ for sufficiently large $m\in \N$.

In case (b) we obtain
\begin{eqnarray*}
\lefteqn{-\langle \left(\Lambda'(\sigma)\kappa \right) g_m,g_m\rangle=\int_\Omega \kappa |\nabla u_\sigma^{g_m}|^2 \dx}\\
&= \int_{D_1} \kappa |\nabla u_{\sigma}^{g_m}|^2 \dx  
+ \int_{D_2} \kappa |\nabla u_{\sigma}^{g_m}|^2 \dx 
\int_{\Omega\setminus (D_1\cup D_2)} \kappa |\nabla u_{\sigma}^{g_m}|^2 \dx\\
&\leq -\mathrm{ess\,inf}\, (-\kappa)|_{D_1} \int_{D_1} |\nabla u_{\sigma}^{g_m}|^2 \dx
+ \norm{\kappa}_{L^\infty(D_2)} \int_{D_2} |\nabla u_{\sigma}^{g_m}|^2 \dx\to -\infty,
\end{eqnarray*}
so that $\int_\Omega \kappa |\nabla u_\sigma^{g_m}|^2 \dx<0$ for sufficiently large $m\in \N$.
\hfill $\Box$

\begin{remark}
It is known (see, e.g., \cite[Cor.~3.5(b)]{harrach2010exact}) that for all piecewise analytic $\sigma$, the Fr\'echet derivative $\Lambda'(\sigma)$ is injective on the space of piecewise analytic functions, i.e. 
$\Lambda'(\sigma)\kappa\neq 0$ for all piecewise analytic $0\not\equiv \kappa\in L^\infty(\Omega)$.

Since $\Lambda'(\sigma)\kappa$ is a compact self-adjoint operator, this means that $\Lambda'(\sigma)\kappa$ must possess 
either a positive or a negative eigenvalue. Lemma \ref{lemma:linearized_uniform_sign} can be interpreted in the
sense, that for each $\kappa\not\equiv 0$ this property is sign-uniform in $\sigma$, i.e., 
for each $\kappa\not\equiv 0$, the operator $\Lambda'(\sigma)\kappa$ either possesses a positive eigenvalue for all $\sigma$,
or it possesses a negative eigenvalue for all $\sigma$ (or both properties are fulfilled). 
\end{remark}

With these preparations we can now show the Theorems \ref{thm:stability_full} and \ref{thm:stability_Galerkin}.

\paragraph{Proof of Theorem \ref{thm:stability_full}.}
The assertion follows from lemma \ref{lemma:estimate_with_f}--\ref{lemma:linearized_uniform_sign} with
\[
c:=\sup_{g\in \Ld^2(\Sigma), \norm{g}=1} f(\hat \tau_1,\hat \tau_2,\hat \kappa,g)>0.
\]
\hfill $\Box$

\paragraph{Proof of Theorem \ref{thm:stability_Galerkin}.}
Using that 
\[
\norm{P_{G_n}\left(\Lambda(\sigma_1)-\Lambda(\sigma_2)\right) P_{G_n}^*}= 
=\sup_{g\in G_n, \norm{g}=1} \left| \Langle g, \left(\Lambda(\sigma_1)-\Lambda(\sigma_2)\right)g \Rangle \right|,
\]
we obtain as in in lemma \ref{lemma:estimate_with_f} and lemma \ref{CM_inf_is_min} that
for all $n\in \N$, there exists $(\hat \tau_1^{(n)},\hat \tau_2^{(n)},\hat \kappa^{(n)}) \in \FF_{[a,b]}\times \FF_{[a,b]}\times \KK$ so that
\begin{equation}\labeq{Galerkin_hilf1}
\frac{\norm{P_{G_n}\left(\Lambda(\sigma_1)-\Lambda(\sigma_2)\right) P_{G_n}^*}}
{\norm{\sigma_1-\sigma_2}}
\geq  \sup_{g\in G_n, \norm{g}=1} f(\hat \tau_1^{(n)},\hat \tau_2^{(n)},\hat \kappa^{(n)},g).
\end{equation}
The right hand side of \req{Galerkin_hilf1} is monotonically increasing in $n\in \N$ since the spaces $G_n$ are nested. 
Hence, the assertion of Theorem \ref{thm:stability_Galerkin} follows, if we can prove that there exists $n\in \N$ with
\begin{equation}\labeq{sup_larger_zero_Galerkin}
\sup_{g\in G_n, \norm{g}=1} f( \tau_1, \tau_2, \kappa,g)>0 \quad \mbox{ for all } 
(\tau_1,\tau_2,\kappa)\in \FF_{[a,b]}\times \FF_{[a,b]}\times \KK.
\end{equation}
We argue by contradiction and assume that this is not the case. Then there exists a sequence 
$(\tau_1^{(n)},\tau_2^{(n)},\kappa^{(n)})\in \mathcal F_{[a,b]}\times \mathcal F_{[a,b]}\times \KK$
with 
\[
\sup_{g\in G_n, \norm{g}=1} f(\tau_1^{(n)},\tau_2^{(n)},\kappa^{(n)},g)\leq 0 \quad \mbox{ for all } n\in \N,
\]
which also implies
\[
\sup_{g\in G_m, \norm{g}=1} f(\tau_1^{(n)},\tau_2^{(n)},\kappa^{(n)},g)\leq 0 \quad \mbox{ for all } n\in \N,\ n\geq m.
\]
After passing to a subsequence if necessary, we can assume by compactness that 
the sequence $(\tau_1^{(n)},\tau_2^{(n)},\kappa^{(n)})$ converges against some element 
\[
(\hat\tau_1,\hat\tau_2,\hat\kappa)\in \FF_{[a,b]}\times \FF_{[a,b]}\times \KK.
\]
Since, for all $m\in \N$, the function
\[
(\tau_1,\tau_2,\kappa)\mapsto \sup_{g\in G_m, \norm{g}=1} f(\tau_1,\tau_2,\kappa,g)
\]
is lower semicontinuous, it follows that
\[
\sup_{g\in G_m, \norm{g}=1} f(\hat\tau_1,\hat\tau_2,\hat\kappa,g)\leq 0 \quad \mbox{ for all } m\in \N.
\]
But, by continuity, this would imply
\[
f(\hat\tau_1,\hat\tau_2,\hat\kappa,g)\leq 0 \quad \mbox{ for all } g\in \overline{\bigcup_{m\in \N} G_m}=\Ld^2(\Sigma),
\]
which contradicts lemma~\ref{lemma:linearized_uniform_sign}. This shows that \req{sup_larger_zero_Galerkin} must be true for sufficiently large $n\in \N$ and
thus theorem \ref{thm:stability_Galerkin} is proven.\hfill $\Box$

\section{Uniqueness and Lipschitz stability from electrode measurements}\label{section:CEM}
\subsection{Setting and main results}\label{subsect:setting}

Now we consider the Complete Electrode Model (CEM). As before let $\sigma\in L^\infty_+(\Omega)$ denote
the conductivity distribution in a smoothly bounded domain $\Omega\subseteq \R^d$, $d\geq 2$. 
We assume that $M\in \N$ open, connected, mutually disjoint electrodes $E_m\subseteq \partial \Omega$, $m=1,\ldots,M$, are attached to the boundary of $\Omega$ with all electrodes having the same contact impedance $z>0$. When a current with strength $J_m\in \R$ is driven through the $m$-th electrode
(with $\sum_{m=1}^M J_m=0$), the resulting electrical potential $(u,U)\in H^1(\Omega)\times \R^M$ solves the following equations:
\begin{eqnarray}
\labeq{CEM1} \nabla \cdot \sigma\nabla u=0 \quad & \mbox{in $\Omega$,}\\
\labeq{CEM2} \sigma\partial_\nu u=0 \quad & \mbox{on $\partial \Omega\setminus \bigcup_{m=1}^M E_m$,}\\
\labeq{CEM3} u+z \sigma \partial_\nu u=\mathrm{const.}=:U_m \quad & \mbox{on $E_m$, $m=1,\ldots,M$,}\\
\labeq{CEM4} \int_{E_m} \sigma \partial_\nu u|_{E_m}\dx[s] =J_m \quad & \mbox{on $E_m$, $m=1,\ldots,M$,}
\end{eqnarray}
where $U=(U_1,\ldots,U_M)\in \R^M$ is a vector containing the electric potentials on the electrodes $E_1,\ldots,E_M$.

It can be shown that \req{CEM1}--\req{CEM4} possess a solution $(u,U)\in H^1(\Omega)\times \R^m$ and that the solution is unique 
under the additional gauge (or ground level) condition $U\in \Rd^M$, where $\Rd^M$ is the subspace of vectors in $\R^m$ with zero mean,
cf., e.g., \cite{somersalo1992existence}.
We can thus define the $M$-electrode current-to-potential operator
\[
R_M(\sigma):\ \Rd^M\to \Rd^M:\ I=(I_1,\ldots,I_M) \mapsto U=(U_1,\ldots,U_M),
\]
where $(u,U)\in H^1(\Omega)\times \Rd^m$ solves \req{CEM1}--\req{CEM4}. Note also that \req{CEM1}--\req{CEM4} are equivalent to the 
variational formulation that $(u,U)\in H^1(\Omega)\times \Rd^m$ solves
\begin{equation}\labeq{CEM_Variational}
\int_\Omega \sigma \nabla u \cdot \nabla w \dx + \sum_{m=1}^M \int_{E_m} \frac{1}{z} (u-U_m)(w-W_m) \dx[s] = \sum_{m=1}^M J_m W_m
\end{equation}
for all $(w,W)\in H^1(\Omega)\times \Rd^m$, cf., again, \cite{somersalo1992existence}.

As in the previous section, we will consider conductivities that belong to a finite dimensional subset 
of piecewise-analytic functions. Additionally, in order to use results from \cite{garde2017convergence} on the approximation properties of the CEM, we assume that the background conductivity is an a-priori known smooth function in a fixed neighborhood $U$ of the boundary $\partial \Omega$, i.e., we assume that $\FF$ is a finite dimensional subset of piecewise-analytic functions, so that there exists $\sigma_0\in C^\infty(\overline U)$ with $\sigma|_U=\sigma_0|_U$ for all $\sigma\in \FF$.
Together with the assumption of a-priori known bounds, we assume (for $b>a>0$)
\[
\sigma\in \FF_{[a,b]}:=\{ \sigma\in \FF:\ a\leq  \sigma(x)\leq  b\mbox{ for all } x\in \Omega\}.
\]
We will show that $R_M(\sigma)$ uniquely determines $\sigma\in\FF_{[a,b]}$ (with Lipschitz stability) if, roughly speaking, enough electrodes are being used. 
To make this statement precise, assume that the number of electrodes is increased so that the electrode configurations fulfill the Hyv\"onen criteria \cite{hyvonen2009approximating,garde2017convergence}:
\begin{description}
\item[(H1)] For all electrode configurations 
\[
E_1,\ldots,E_M\subseteq \partial \Omega
\]
there exist open, connected, and mutually disjoint sets (called virtual extended electrodes) $\tilde E_m^{(M)}$, $m=1,\ldots,M$ with
\begin{eqnarray*}
E_m^{(M)}\subseteq \tilde E_m^{(M)}, \quad \bigcup_{m=1}^M \overline{\tilde E_m^{(M)}}&=\partial \Omega,
\end{eqnarray*}
so that
\begin{eqnarray*}
h_{M}:=\max_{m=1,\ldots,M} \left\{\sup_{x,y\in \tilde E_m^{(M)}} \mathrm{dist}(x-y)\right\} \to 0 \quad & \mbox{ for } M\to \infty,\\
\exists c_{\mathrm{E}}>0:\ \min_{m=1,\ldots,M} \frac{|E_m^{(M)}|}{|\tilde E_m^{(M)}|}\geq c_{\mathrm{E}} \quad & \mbox{ for all } M\in \N.
\end{eqnarray*}
\item[(H2)] The operators
\begin{eqnarray*}
Q_{M}:& \ \R^M\to L^2(\partial \Omega),\quad (J_m)_{m=1}^M \mapsto \sum_{m=1}^M J_m \chi_{\tilde E_m^{(M)}}\\
P_{M}:& \ L^2(\partial \Omega)\to \R^M,\quad g\mapsto \left(\frac{1}{|E_m^{(M)}|} \int_{E_m^{(M)}} g \dx[s]\right)_{m=1}^M
\end{eqnarray*}
fulfill that
\begin{eqnarray*}
\exists C_{\mathrm{E}}>0:\ & \norm{(I-Q_{M}P_{M})f}_{L^2(\partial \Omega)}\leq C_{\mathrm{E}} h_M \inf_{c\in \R} \norm{f+c}_{H^1(\partial \Omega)}
\end{eqnarray*}
for all $M\in \N$ and all $f\in H^1(\partial \Omega)$.
\end{description}
The first criterion implies the natural assumption that the electrode sizes shrink to zero, but always cover a certain fraction of the boundary. 
The somewhat technical second criterion can be interpreted as a Poincaré-type inequality that is fulfilled for regular enough electrode shapes, see \cite{lechleiter2008factorization,hyvonen2009approximating,garde2017convergence}. Together these criteria guarantee that the electrode measurement approximate all possible continuous measurements in a suitable sense.

Now we can state our main result:
\begin{theorem}\label{thm:stability_CEM}
There exists $N\in \N$ and $c>0$ such that for all $M\geq N$
\[
\norm{R_M(\sigma_1) - R_M(\sigma_2)}_{\LL(\Rd^M)}\geq c \norm{\sigma_1-\sigma_2}_{L^\infty(\Omega)} \quad \mbox{ for all } \sigma_1,\sigma_2\in \FF_{[a,b]}
\]
In particular, this implies that for all $\sigma_1,\sigma_2\in \FF_{[a,b]}$ and $M\geq N$
\[
R_M(\sigma_1) = R_M(\sigma_2)  \quad \mbox{ if and only if } \quad \sigma_1=\sigma_2.
\]
\end{theorem}
Theorem \ref{thm:stability_CEM} will be proven in subsection~\ref{subsect:proofs_CEM}.

\subsection{Differentiability, monotonicity, and approximation of linearized measurements}\label{subsect:CEM_tools}
The electrode measurements $R_M(\sigma)$ fulfill analogue differentiability and monotonicity properties as the Neumann-to-Dirichlet-Operators. 
In the following $\langle \cdot,\cdot \rangle_M$ denotes the Euclidian scalar product in $\R^M$. For a vector $J=(J_1,\ldots,J_M)\in \Rd^M$, we denote by $u^{(J)}_\sigma\in H^1_\diamond(\Omega)$ the solution of the CEM equations \req{CEM1}--\req{CEM4} with conductivity $\sigma\in L^\infty_+(\Omega)$ and electrode currents $J_1,\ldots,J_m\in \R$.

\begin{lemma}\label{lemma:Frechet_differentiable_CEM}
\begin{enumerate}[(a)]
\item The mapping 
\[
R_M:\ L^\infty_+(\Omega)\to \LL(\Rd^M),\quad \sigma\mapsto R_M(\sigma)
\]
is Fr\'echet differentiable. Its derivative is given by
\[
R_M'(\sigma)\in \LL(L^\infty(\Omega), \LL(\Rd^M)), \quad 
\left(R_M'(\sigma)\kappa\right) J= V,
\]
where $(v,V)\in H^1(\Omega)\times \Rd^m$ solves
\begin{eqnarray}
\nonumber \int_\Omega \sigma \nabla v \cdot \nabla w \dx + \sum_{m=1}^M \int_{E_m} \frac{1}{z} (v-V_m)(w-W_m) \dx[s]\\
 = - \int_\Omega \kappa \nabla u_\sigma^{(J)}\cdot
\nabla w \dx \labeq{Diff_CEM}
\end{eqnarray}
for all $(w,W)\in H^1(\Omega)\times \Rd^m$.
\item For all $\sigma\in L^\infty_+(\Omega)$ and $\kappa\in L^\infty(\Omega)$ the operator $R_M'(\sigma)\kappa\in  \LL(\Rd^M)$ is self-adjoint,
and it fulfills
\[
\langle \left(R_M'(\sigma)\kappa\right) I,J\rangle_M=-\int_\Omega \kappa \nabla u_\sigma^{(I)} \cdot \nabla u_\sigma^{(J)} \dx.
\]
for all $I,J\in \Rd^M$.
\item The mapping
\[
R_M':\ L^\infty_+(\Omega)\to \LL(L^\infty(\Omega), \LL(\Rd^M)),\quad \sigma\mapsto R_M'(\sigma)
\]
is continuous.
\end{enumerate}
\end{lemma}
\paragraph{Proof.}
This follows from the variational formulation of the CEM \req{CEM_Variational}, cf., e.g., \cite[Section~2]{lechleiter2008newton} or \cite[Appendix~B]{garde2017convergence}.
\hfill $\Box$

\begin{lemma}\label{lemma:CEM_monotony}
For all $\sigma_1,\sigma_2\in L^\infty_+(\Omega)$ and $J\in \Rd^M$, it holds that
\begin{eqnarray*}\label{eq:CEM_monotony_inequality}
 \langle \left(R'(\sigma_2)(\sigma_1-\sigma_2)\right) J,J\rangle_M&=
\int_\Omega(\sigma_2-\sigma_1) |\nabla u_{\sigma_2}^{(J)}|^2 \dx\\ 
&\leq \Langle \left(R_M(\sigma_1)-R_M(\sigma_2)\right)J,J \Rangle_M.
\end{eqnarray*}
\end{lemma}
\paragraph{Proof.}
\cite[Theorem~2]{harrach2015resolution}.
\hfill $\Box$

We will also require the following result from Garde and Staboulis \cite{garde2017convergence}
that the linearized CEM measurements approximate the linearized Neumann-to-Dirichlet operator.
\begin{lemma}\label{lemma:approx_linearization}
Under the Hyv\"onen assumptions (H1) and (H2), there exists $C>0$ such that for all $\sigma\in \FF_{[a,b]}$ and $\kappa\in \span \FF$
\[
\norm{\Lambda'(\sigma)\kappa- L Q \left( R'(\sigma)\kappa\right) Q^* }_{\LL(\Ld^2(\partial \Omega))}
\leq C h_M \norm{\sigma}_{L^\infty(\Omega)} \norm{\kappa}_{L^\infty(\Omega)},
\]
where
\begin{eqnarray*}
L:\   L^2(\partial \Omega)&\to L^2_\diamond(\partial \Omega),\quad & g \mapsto g-\frac{1}{|\partial \Omega|}\int_{\partial \Omega} g \dx[s],\\
Q_M^*:\ L^2(\partial \Omega)&\to \R^M,\quad & g \mapsto   \left(\int_{\tilde E_1^{(M)}} g\dx[s],\ldots,\int_{\tilde E_M^{(M)}} g\dx[s]\right).
\end{eqnarray*}
Moreover, for all $g\in \Ld^2(\partial \Omega)$
\[
\Langle L Q \left( R'(\sigma)\kappa\right) Q^* g, g\Rangle = \Langle \left( R'(\sigma)\kappa\right) Q^* g, Q^* g\Rangle_M.
\]
\end{lemma}
\paragraph{Proof.}
\cite[Thm.~3,\ Prop.~4]{garde2017convergence}
\hfill $\Box$

\subsection{Proof of Theorem \ref{thm:stability_CEM}}\label{subsect:proofs_CEM}

Again, for the sake of brevity, we omit norm subscripts when the choice of the norm is clear from the context.
As in subsection \ref{subsect:proofs} we obtain from the monotonicity result \ref{lemma:CEM_monotony} that
\begin{eqnarray*}
\norm{R_M(\sigma_1)-R_M(\sigma_2)} \geq  \norm{\sigma_1-\sigma_2} \inf_{(\tau_1,\tau_2,\kappa) \atop \in \FF_{[a,b]}\times \FF_{[a,b]}\times \KK} \ \sup_{J\in \Rd^M\atop \norm{J}=1} f_M(\tau_1,\tau_2,\kappa,J).
\end{eqnarray*}
where $f_M:\ L^\infty_+(\Omega)\times L^\infty_+(\Omega)\times L^\infty(\Omega)\times \Rd^M\to \R$ is defined by
\[
f_M(\tau_1,\tau_2,\kappa,J):=\max \left\{ \Langle \left(R_M'(\tau_1) \kappa\right) J,J \Rangle,
-\Langle \left(R_M'(\tau_2) \kappa\right) J,J \Rangle\right\},
\]
and $\KK:=\{ \kappa\in \span \FF:\ \norm{\kappa}=1\}$.

We compare this with $f:\ L^\infty_+(\Omega)\times L^\infty_+(\Omega)\times L^\infty(\Omega)\times \Ld^2(\partial \Omega)\to \R$,
\[
f(\tau_1,\tau_2,\kappa,g):=\max \left\{ \Langle \left(\Lambda'(\tau_1) \kappa\right) g,g \Rangle,
-\Langle \left(\Lambda'(\tau_2) \kappa\right) g,g \Rangle\right\},
\]
from the continuum model (cf.\ Lemma~\ref{lemma:estimate_with_f}) with $\Sigma=\partial \Omega$. We obtain with lemma~\ref{lemma:approx_linearization} that for all $(\tau_1,\tau_2,\kappa)  \in \FF_{[a,b]}\times \FF_{[a,b]}\times \KK$ and $g\in \Ld^2(\partial \Omega)$
with $\norm{g}=1$.
\begin{eqnarray*}
\left| f(\tau_1,\tau_2,\kappa,g)-f_M(\tau_1,\tau_2,\kappa,Q_M^* g)\right|
%
& \leq Ch_M \norm{\kappa} \norm{g}
\max\{ \norm{\tau_1},\norm{\tau_2}\}\\
&\leq C h_M b.
\end{eqnarray*}
For all $g\in \R^M$ we have that
\begin{eqnarray*}
\norm{Q_M ^* g}^2&= \sum_{m=1}^M  \left(\int_{\tilde E_m^{(M)}} g \dx[s]\right)^2
\leq \sum_{m=1}^M  |\tilde E_m^{(M)}| \int_{\tilde E_m^{(M)}} g^2 \dx[s]
\leq \left| \partial \Omega \right| \norm{g}^2,
\end{eqnarray*}
so that we obtain for all $(\tau_1,\tau_2,\kappa)  \in \FF_{[a,b]}\times \FF_{[a,b]}\times \KK$
\begin{eqnarray*}
\sup_{J\in \Rd^M,\atop \norm{J}=1} f_M(\tau_1,\tau_2,\kappa,J)
& \geq \left| \partial \Omega \right|^{-1} \sup_{g\in \Ld^2(\partial \Omega),\atop \norm{g}=1} f_M(\tau_1,\tau_2,\kappa,Q_M^*g)\\
& \geq \left| \partial \Omega \right|^{-1} \sup_{g\in \Ld^2(\partial \Omega),\atop \norm{g}=1} f(\tau_1,\tau_2,\kappa,g) - \left| \partial \Omega \right|^{-1} Ch_M b.
\end{eqnarray*}
Since the first summand is positive by lemma~\ref{lemma:linearized_uniform_sign}, it follows that
for sufficiently large numbers of electrodes $M$ 
\[
\sup_{J\in \Rd^M,\atop \norm{J}=1} f_M(\tau_1,\tau_2,\kappa,J)>0 \quad \mbox{ for all } (\tau_1,\tau_2,\kappa)  \in \FF_{[a,b]}\times \FF_{[a,b]}\times \KK.
\]
With the same lower semicontinuity and compactness argument as in the continuum model, this yields
\[
\inf_{(\tau_1,\tau_2,\kappa) \atop \in \FF_{[a,b]}\times \FF_{[a,b]}\times \KK} \ \sup_{J\in \Rd^M,\atop \norm{J}=1} f_M(\tau_1,\tau_2,\kappa,J)>0,
\]
so that the assertion is proven.\hfill $\Box$

\ack

This work is devoted to the memory of Professor Armin Lechleiter who will be deeply missed as a scientist, and as a friend.

\section*{References}

\bibliographystyle{abbrv}
\bibliography{literaturliste}

\end{document}